\documentclass{elsart}
\usepackage{amsfonts,amsmath,natbib}
\usepackage{graphicx}

\newcommand{\Tr}{\operatorname{Tr}}
\newcommand{\hf}{\hat{f}}

\newcommand{\eqn}[1]{(\ref{#1})}

\newcommand{\cc}{{\mathbb C}}

\newcommand{\bb}{{\bf b}}

\newcommand{\tabl}[1]{{\left[\begin{smallmatrix} #1 \end{smallmatrix}\right]}}
\begin{document}
\begin{frontmatter}

\title{Markov Bases for Noncommutative Fourier Analysis of Ranked Data}
\author[D]{Persi Diaconis},
\ead{diaconis@math.stanford.edu}
\author[E]{Nicholas Eriksson}%\corauthref{cor}}
%\corauth[cor]{Corresponding author.}
\ead{eriksson@math.berkeley.edu}

\date{October 15, 2004}

\address[D]{Departments of Mathematics and Statistics, 
Stanford University}%, Palo Alto, CA}%94305-4065}
\address[E]{Department of Mathematics, 
University of California at Berkeley}% 94720-3840}

\begin{abstract}
To calibrate Fourier analysis of $S_5$ ranking data by Markov chain
Monte Carlo techniques, a set of moves (Markov basis) is needed.  
We calculate this basis, and use it to provide a new
statistical analysis of two data sets.  The calculation
involves a large Gr\"obner basis computation (45825 generators), but
reduction to a minimal basis and reduction by natural symmetries leads to a remarkably small
basis (14 elements).  Although the Gr\"obner basis calculation is infeasible for $S_6$, 
we exploit the symmetry of the problem to calculate a Markov basis for $S_6$ with  
7,113,390 elements in 58 symmetry classes.
We improve a bound on the degree of the generators for a
Markov basis for $S_n$ and conjecture that this ideal is
generated in degree 3.

\end{abstract}

%\begin{keyword}
%Markov chain Monte Carlo,
%Markov basis, ranked data, non-commutative Fourier analysis, toric ideals.
%\MSC MSC code \sep MSC code2
%6E17 = ???
%13P10 = Polynomial ideals, Gr¨obner bases

%\end{keyword}

\end{frontmatter}

\section{Election data with five candidates}
Table~\ref{tab:1} shows the results of an election.  A population of 5738
voters was asked to rank five candidates for president of a national
professional organization.  The table shows the number of voters
choosing each ranking.  
\begin{table}
\caption{American Psychological Association voting data: the number of voters that ranked the 5 candidates in a given order.}
\label{tab:1}
\centering
\begin{tabular}{|cc|cc|cc|cc|}
\hline
Ranking & \# votes & Ranking & \# votes & Ranking & \# votes & Ranking & \# votes\\ \hline
54321 & 29 & 43521 & 91 & 32541 & 41 & 21543 & 36 \\
54312 & 67 & 43512 & 84 & 32514 & 64 & 21534 & 42 \\
54231 & 37 & 43251 & 30 & 32451 & 34 & 21453 & 24 \\
54213 & 24 & 43215 & 35 & 32415 & 75 & 21435 & 26 \\
54132 & 43 & 43152 & 38 & 32154 & 82 & 21354 & 30 \\
54123 & 28 & 43125 & 35 & 32145 & 74 & 21345 & 40 \\
53421 & 57 & 42531 & 58 & 31542 & 30 & 15432 & 40 \\
53412 & 49 & 42513 & 66 & 31524 & 34 & 15423 & 35 \\
53241 & 22 & 42351 & 24 & 31452 & 40 & 15342 & 36 \\
53214 & 22 & 42315 & 51 & 31425 & 42 & 15324 & 17 \\
53142 & 34 & 42153 & 52 & 31254 & 30 & 15243 & 70 \\
53124 & 26 & 42135 & 40 & 31245 & 34 & 15234 & 50 \\
52431 & 54 & 41532 & 50 & 25431 & 35 & 14532 & 52 \\
52413 & 44 & 41523 & 45 & 25413 & 34 & 14523 & 48 \\
52341 & 26 & 41352 & 31 & 25341 & 40 & 14352 & 51 \\
52314 & 24 & 41325 & 23 & 25314 & 21 & 14325 & 24 \\
52143 & 35 & 41253 & 22 & 25143 & 106 & 14253 & 70 \\
52134 & 50 & 41235 & 16 & 25134 & 79 & 14235 & 45 \\
51432 & 50 & 35421 & 71 & 24531 & 63 & 13542 & 35 \\
51423 & 46 & 35412 & 61 & 24513 & 53 & 13524 & 28 \\
51342 & 25 & 35241 & 41 & 24351 & 44 & 13452 & 37 \\
51324 & 19 & 35214 & 27 & 24315 & 28 & 13425 & 35 \\
51243 & 11 & 35142 & 45 & 24153 & 162 & 13254 & 95 \\
51234 & 29 & 35124 & 36 & 24135 & 96 & 13245 & 102 \\
45321 & 31 & 34521 & 107 & 23541 & 45 & 12543 & 34 \\
45312 & 54 & 34512 & 133 & 23514 & 52 & 12534 & 35 \\
45231 & 34 & 34251 & 62 & 23451 & 53 & 12453 & 29 \\
45213 & 24 & 34215 & 28 & 23415 & 52 & 12435 & 27 \\
45132 & 38 & 34152 & 87 & 23154 & 186 & 12354 & 28 \\
45123 & 30 & 34125 & 35 & 23145 & 172 & 12345 & 30 \\\hline
\end{tabular}
\end{table}
For example, 29 voters ranked candidate 5 first, candidate 4 second,
\dots, and candidate 1 last, resulting in the entry $54321=29$. 
Table~\ref{tab:firstorder} shows a simple summary of the data:
the proportion of voters ranking candidate $i$ in position $j$.  For
example, 28.0\% of the voters ranked candidate 3 first and 23.1\% of
the voters ranked candidate 3 last.

\begin{table}
\caption{First order statistics: The proportion of voters who ranked candidate $i$ in position $j$.
This is a scaled version of the Fourier transform of Table~\ref{tab:1} at the permutation representation.}
\label{tab:firstorder}
\centering
\begin{tabular}{|c|ccccc|}
\hline
&&Rank & & &\\
Candidate & 1&2&3&4&5\\\hline
1&18.3& 26.4& 22.8& 17.4& 14.8\\
2&13.5& 18.7& 24.6& 24.6& 18.3\\
3&28.0& 16.7& 13.8& 18.2& 23.1\\
4&20.4& 16.9& 18.9& 20.2& 23.3\\
5&19.6& 21.0& 19.6& 19.2& 20.3\\\hline
\end{tabular}
\end{table}

Table~\ref{tab:firstorder} is a natural summary of the 120 numbers in
Table~\ref{tab:1}, but is it an adequate summary?  Does it capture all
the ``juice'' in the data?  In this paper, we develop tools to answer
such questions using Fourier analysis and algebraic techniques.

In Section~\ref{sec:fourier}, we give a general exposition of how
noncommutative Fourier analysis can be used to analyze group valued
data with summary given by a representation $\rho$.  In order to use
Markov chain Monte Carlo techniques to calibrate the Fourier analysis,
we define an exponential family and toric ideal associated to a finite
group $G$ and integer representation $\rho$.  A generating set of the
toric ideal can be used to run a Markov chain to sample from data on
the group.  For example, the 14 moves in Table~\ref{tab:s5moves} allow
us to randomly sample from the space of data on $S_5$ with fixed first
order summary (Table~\ref{tab:firstorder}).

In Section~\ref{sec:compute} we show how this basis (Table~\ref{tab:s5moves}) was computed -
either using Gr\"obner bases or by utilizing symmetry.  We
describe extensive computations of the basis for ranked data on at most
6 objects.  From these computations, we conjecture that the toric ideal
for $S_n$  is generated in degree 3.
In Section~\ref{sec:math}, we show this ideal for $S_n$ is  generated in degree $n-1$, improving
a result of \cite{ds}, and we describe the degree 2 moves.
Finally, in Section~\ref{sec:stat}, we apply these methods to analyze the data in Table~\ref{tab:1} 
and an example from \cite{ds}.

\begin{table}
\caption{$S_5$ moves: there are 29890 moves in 14 symmetry classes of sizes 200-7200}
\label{tab:s5moves}
\centering
\begin{tabular}{|cc|cc|}
\hline
Move & Number & Move & Number \\\hline
$\tabl{53412\\54321}-\tabl{53421\\54312}$ & 450&
$\tabl{45231\\54312}-\tabl{45312\\54231}$ & 600\\
$\tabl{54123\\54231\\54312}-\tabl{54132\\54213\\54321}$ & 200&
$\tabl{53412\\54123\\54231}-\tabl{53421\\54132\\54213}$ & 3600\\
$\tabl{45123\\54231\\54312}-\tabl{45132\\54213\\54321}$ & 200&
$\tabl{45123\\53412\\54231}-\tabl{45132\\53421\\54213}$ & 7200\\
$\tabl{43512\\54123\\54231}-\tabl{43521\\54132\\54213}$ & 3600&
$\tabl{43512\\53241\\54123}-\tabl{43521\\53142\\54213}$ & 3600\\
$\tabl{45231\\52341\\53412}-\tabl{45312\\52431\\53241}$ & 7200&
$\tabl{45132\\52341\\53412}-\tabl{45312\\52431\\53142}$ &3600\\
$\tabl{34512\\45123\\53241}-\tabl{34521\\45213\\53142}$ &600&
$\tabl{34521\\45213\\53142}-\tabl{35142\\43521\\54213}$ &600\\
$\tabl{35142\\43521\\54213}-\tabl{35241\\43512\\54123}$ &600&
$\tabl{34521\\45312\\52143}-\tabl{35142\\42513\\54321}$ & 1440\\
\hline
\end{tabular}
\end{table}

\section{Fourier analysis of group valued data} \label{sec:fourier}
Let $G$ be a finite group (in our example, $G = S_5$).  Let $f \colon
G \to \Zset$ be any function on $G$.  For example, if $g_1, g_2,
\dots, g_N$ is a sample of points chosen from a distribution on G,
take $f(g)$ to be the number of sample points $g_i$ that are equal to $g$.
We view $f$ interchangeably as either a function on the group or an
element of the group ring $\Zset[G]$.
Recall that a map $\rho \colon G \to GL(V_\rho)$ is a matrix
representation of $G$ if $\rho(st) = \rho(s)\rho(t)$ for all $s,t \in
G$.  The dimension $d_\rho$ of the representation $\rho$ is the
dimension of $V_\rho$ as a $\cc$-vector space.
We say that a  $\rho$ is integer-valued if 
$\rho_{ij}(g) \in \Zset$ for all $g\in G$ and for all $1 \leq i,j \leq d_\rho$.
We denote the set of irreducible representations of $G$ by $\hat{G}$.

An analysis of $f(g)$ may be based on the Fourier transform.  
The Fourier transform of $f$ at $\rho$ is
\begin{equation}
\hat{f}(\rho) = \sum_{g \in G} f(g) \rho(g).
\end{equation}
The Fourier transform at all the irreducible representations $\rho \in
\hat{G}$ determines $f$ through the Fourier inversion formula
\begin{equation}\label{eq:inv}
%also f(g) =  \sum_{\chi} \chi(1) \chi(f g^{-1})
f(g) = \frac 1 {|G|} 
\sum_{\rho \in \hat{G}} d_{\rho} \Tr(\hat{f}(\rho) \rho(g^{-1})),
\end{equation}
which can be rewritten as
$f(g)= \sum_{\rho\in \hat{G}} f|_{V_\rho}(g)$,
where
\begin{equation} \label{eq:proj}
f|_{V_\rho}(g) = 
\frac {d_\rho} {|G|} \sum_{h \in G} \chi_\rho(h) f(g h).
\end{equation}
This decomposition shows the contributions to $f$ from 
each of the irreducible representations of $G$.
For example, if a few of the $f|_{V_\rho}$ are large, we can analyze these
components in order to understand the structure of $f$.  See
\cite{diaconis-book,diaconis-spectral} for background, proofs, and
previous literature.

\begin{exmp}
This analysis 
is most familiar for the cyclic group $C_n$ where
it becomes the discrete Fourier transform
\begin{equation}\label{eq:classical}
\hf(j) = \sum_{k=0}^{n-1}f(k) e^{-2 \pi i j k / n}, \qquad
f(k) = \frac 1 n \sum_{j=0}^{n-1} \hf(j) e^{2 \pi i k j / n} 
\end{equation}
In \eqn{eq:classical}, if a few of the $\hf(j)$ are much larger than
the rest, then $f$ is well understood as approximately a sum of a few
periodic components. 
\end{exmp}

For the symmetric group $S_n$, the  permutation representation
assigns permutation matrices $\rho(\pi)$ to permutations $\pi$.  Thus, if
$f(\pi)$ is the number of rankers choosing $\pi$, $\hf(\rho)$ is a $n
\times n$ matrix with $(i,j)$ entry the number of rankers ranking item
$i$ in position $j$ (as in Table~\ref{tab:firstorder}).  
The irreducible representations of $S_5$ are indexed by the seven
partitions of five and are written as $S^{\lambda}$ where $\lambda$ is a partition of 5.  
For our data,  \eqn{eq:inv} gives a decomposition of $f$ into 7 parts.
Table~\ref{tab:sum} shows the lengths of the projection of
Table~\ref{tab:1} onto the seven isotypic subspaces of $S_5$.

\begin{table}
\caption{Squared length (divided by 120) of the projection of the APA data into the 7 isotypic subspaces of $S_5$.}
\label{tab:sum}
\centering
     \begin{tabular}{|c|ccccccc|}
	\hline
	& $S^{5}$ & $S^{4,1}$ &  $S^{3,2}$ & $S^{3,1,1}$ & 
	$S^{2,2,1}$ & $S^{2,1,1,1}$ & $S^{1,1,1,1,1}$\\
	$d_\rho^2$  & 1 & 16 & 25 & 36 & 25 & 16 & 1\\ \hline
	Data      & 2286 & 298 & 459 & 78 & 27 & 7 & 0\\\hline
      \end{tabular}
\end{table}

The largest contribution to 
the data occurs from the trivial representation $S^{5}$.   We call the projection onto
$S^{5} \oplus S^{4,1}$ the first order summary; it was shown in
Table~\ref{tab:firstorder} above.  We see that the projection onto
$S^{3,2}$ is also sizable while the rest of the projections are relatively
negligible.  This suggests a data-analytic look at the projection into
$S^{3,2}$.  
Table~\ref{tab:secondorder} shows this projection
in a natural coordinate system. This projection is based on the permutation
representation of $S_5$ on unordered pairs $\{i,j\}$.  
Table~\ref{tab:secondorder} is an embedding of a 25 dimensional space
into a 100 dimensional space so that its coordinates are easy to
interpret.  See \cite{diaconis-spectral} for further explanation.

\begin{table}
\caption{Second order summary for the APA data}
\label{tab:secondorder}
\centering
\begin{tabular}{|c|cccccccccc|}
\hline
& Rank&&&&&&&&& \\
Candidate & 1,2 & 1,3 & 1,4 & 1,5 & 2,3 & 2,4 & 2,5 & 3,4 & 3,5 & 4,5 \\ \hline
1,2 & -137 & -20 & 18 & 140 & 111 & 22 & 4 & 6 & -97 & -46\\
1,3 & 476  & -88 & -179 & -209 & -147 & -169 & -160 & 107 & 128 & 241\\
1,4 & -189 & 51  & 113 & 24 & -9 & 98 & 99 & -65 & 23 & -146\\
1,5 & -150 & 57  & 47 & 45 & 43 & 49 & 56 & -48 & -53 & -48\\
2,3 & -42  & 84  & 19 & -61 & 30 & -16 & 82 & -76 & -39 & 72\\
2,4 & 157  & -20 & -43 & -25 & -93 & -76 & -56 & 8 & 38 & 112\\
2,5 & 22   & -44 & 7 & 15 & -117 & 69 & 25 & 62 & 99 & -138\\
3,4 & -265 & -7  & 72 & 199 & 39 & 140 & 85 & 19 & -52 & -233\\
3,5 & -169 & 10  & 88 & 70 & 78 & 44 & 47 & -51 & -36 & -80\\
4,5 &296  & -24 & -142 & -130 & -5 & -163 & -128 & 38 & -9 & 267\\
\hline
\end{tabular}
\end{table}

The largest number in Table~\ref{tab:secondorder} is 476 in the
$\{1,3\}, \{1,2\}$ position corresponding to a large positive
contribution to ranking candidates one and three in the top two
positions.  There is also a large positive contribution for ranking
candidates four and five in the top two positions.  Since
Table~\ref{tab:secondorder} gives the projection of $f$ onto a subspace
orthogonal to $S^{5} \oplus S^{4,1}$, the popularity of individual
candidates has been subtracted out.  We can see the ``hate vote''
against the pair of candidates one and three (and the pair four and
five) from the last column.  Finally, the negative entries for e.g.,
pairs one and four, one and five, three and four, three and five show
that voters don't rank these pairs in the same way.

The preceding analysis is from \citet{diaconis-spectral} which used
it to show that noncommutative spectral analysis could be a useful
adjunct to other statistical techniques for data analysis.

The data is from the American Psychological Association -- a polarized
group of academicians and clinicians who are on very uneasy terms (the
organization almost split in two just after this election).
Candidates one and three are in one camp, candidates four and five
from the other.  Candidate two seems disliked by both camps.
The winner of the election depends on the method of allocating votes.
For example, the Hare system or plurality voting would elect candidate
three.  However, other widely used voting methods (Borda's sum of
ranks or Coomb's elimination system) elect candidate one.  For details
and further analysis of the data, see \citet{stern}.

To explain the perturbation analysis in Section~\ref{sec:stat}, it is useful to
consider a simple exponential model for group-valued data.

\begin{defn}\label{def:model}
Let $\rho$ be an integer valued representation of a finite group $G$.
Then the exponential family of $G$ and $\rho$ is given by the
family of probability distributions on $G$
\begin{equation}
P_\Theta (g) = Z^{-1} \e^{\Tr(\Theta \rho(g))}  \label{eq:model}
\end{equation}
where the normalizing constant is $Z = \sum_{g \in G} \e^{\Tr(\Theta\rho(g))}$ and
$\Theta$ is a $n \times n$ matrix of parameters to be chosen to fit the data.  
\end{defn}

For example let $G = S_n$ and $\rho$ be the usual permutation representation.  
Then if $\Theta$ is the zero
matrix, $P_\Theta$ is the uniform distribution.  If $\Theta_{1,1}$ is
nonzero and $\Theta_{i,j}$ is zero otherwise, the model $P_\Theta$
corresponds to item one being ranked first with special probability,
the rest ranked randomly.  Such models have been studied by 
\citet{silver,verd,diaconis-spectral}.  See 
\citet{marden} for a book-length treatment of models for permutation data.

From the Darmois-Koopman-Pitman Theorem \citep[e.g.,][Theorem~3.1]{diaconis-freedman}, we deduce
\begin{prop}
The model \eqn{eq:model} has the property that a sufficient statistic
for $\Theta$ based on data $f(\pi)$ is the Fourier transform
$\hf(\rho)$.  Furthermore, \eqn{eq:model} is the unique model characterized by this property.
\end{prop}

\begin{defn}\label{def:ideal}
Given a finite group $G$ and an integer valued representation $\rho$
of dimension $d_\rho$ 
define the toric ideal of $G$ at $\rho$ as $I_{G,\rho} = \ker(\phi_{G,\rho})$, where
\begin{align*}
\phi_{G,\rho} \colon \Cset[x_g \mid g \in G] 
&\longrightarrow \Cset[t_{ij}^{\pm 1} \mid 1 \leq i,j \leq d_\rho]\\
x_g &\longmapsto \prod_{1 \leq i,j \leq d_\rho} t_{ij}^{\rho_{ij}(g)}.
\end{align*}
\end{defn}
This ideal is the vanishing ideal of the exponential family from Definition~\ref{def:model}.
It will be our main object of study in Sections~\ref{sec:compute} and \ref{sec:math}.

As suggested by \cite{fisher}, tests of goodness of
fit of the model \eqn{eq:model} should be based on the conditional
distribution of the data $f$ given the sufficient statistic
$\hf(\rho)$.  By an elementary calculation,
\begin{equation}\label{eq:hyper}
P_\Theta (f | \hf(\rho)) = w^{-1} \prod_{\sigma \in G} \frac 1 {f(\sigma)!},
%\end{equation}
\quad \text{where} \quad
%\[
w = \sum_{\substack{g\in\Zset[G]\\\hat{g}(\rho) = \hat{f}(\rho)}} \prod_{\sigma \in G} \frac 1 {g(\sigma)!}.
\end{equation}
Observe that the conditional distribution in \eqn{eq:hyper} is free of
the unknown parameter $\Theta$.

The original justification for the Fourier decomposition
is model free (non-parametric).  The first order summary in
Table~\ref{tab:firstorder} is a natural object to look at and the
second order summary was analyzed because of a sizable projection to $S^{3,2}$ in Table~\ref{tab:sum}.
It is natural to wonder if the second order
summary is real or just a consequence of finding patterns in any set
of numbers.  To be honest, the APA data is not a sample (those 5,972
who choose to vote are likely to be quite different from the bulk of
the 100,000 or so APA members).  If the first order summary is
accepted ``as is'', the largest probability model for which
$\hat{f}(\rho)$ captures all the structure in the data is the
exponential family \eqn{eq:model}.  It seems natural to use the
conditional distribution of the data given $\hf(\rho)$ as a way of
perturbing things.  The uniform distribution on data with fixed
$\hf(\rho)$ is a much more aggressive perturbation procedure.  Both
are computed and compared in Section~\ref{sec:stat}.

\section{Computing Markov bases for permutation data}\label{sec:compute}

To carry out a test based on Fisher's principles, we use Markov chain Monte Carlo
to draw samples from the distribution \eqn{eq:hyper}.

\begin{defn}\label{def:5}
A {\em Markov basis} for a finite group $G$ and a representation $\rho$ is
a finite subset of ``moves'' $g_1, \dots, g_B\ \in \Zset[G]$ with $\hat{g}_i(\rho) = 0$ such 
that any two elements in $\Nset[G]$
with the same Fourier transform at the representation $\rho$
can be connected by a sequence of moves in that subset.
\end{defn}

In \cite{ds} it was explained how Gr\"obner basis techniques could be
applied to find such Markov bases.
\begin{prop}
A generating set of $I_{G,\rho}$ (see Definition~\ref{def:ideal}) is a Markov basis
for the group $G$ and the representation $\rho$.
\end{prop}

We will write $I_{S_n}$ for our main example, the ideal of $S_n$ with
the permutation representation $\rho$.
The representation $\rho \colon \Nset[S_n] \to \Nset^{n^2}$ sends an element of $S_n$ to its permutation matrix.
%and extends by linearity.
The elements $\bb \in \Nset^{n^2}$ 
with $\rho^{-1}(\bb)$ non-empty are
the {\em magic squares}, that is, matrices with non-negative integer entries
such that all row and column sum are equal.
We write an element $\pi_1 + \dots + \pi_m  \in \Nset[S_n]$ as a tableau
$\tabl{\pi_1(1) & \dots &\pi_1(n) \\\vdots &  & \vdots \\\pi_m(1) & \dots &\pi_m(n)}$.
In this notation, a Markov basis element is written as a difference of two tableaux.
For example, the degree 2 element of the Markov basis for $S_5$,
$\tabl{13452\\14325} - \tabl{13425\\14352}$,
corresponds to adding one to the entries 13452 and 14325 in
Table~\ref{tab:1} and subtracting one from the entries 13425 and
14352.  

At the time of writing \cite{ds}, finding a
Gr\"obner basis for $I_{S_5}$ was computationally infeasible.
Due to an increase in computing power
and the development of the software {\tt 4ti2} \citep{4ti2}, we were able to
compute a Gr\"obner and a minimal basis of $I_{S_5}$.

This computation involved finding a Gr\"obner basis of a
toric ideal involving 120 indeterminates.  It  took {\tt 4ti2} 
approximately 90 hours of CPU time on a 2GHz machine and 
produced a basis with 45,825 elements.   
The Markov basis  had
29890 elements, 1050 of degree 2 and 28840 of degree 3, see
Tables~\ref{tab:s5moves} and \ref{tab:basis}.
Using {\tt 4ti2}, 
we have also computed Markov bases of the ideals $I_{S_n}$ for $n =
3$ and $n=4$, they are shown in Table~\ref{tab:s4moves}.

\begin{table}
\caption{Markov bases for $S_3$ and $S_4$ and the size of their symmetry classes.}
\label{tab:s4moves}
\centering
\begin{tabular}{|cc|cc|}
\hline
%\multicolumn{2}{c}{$S_3$} & \multicolumn{2}{c}{$S_4$} \\
$S_3$ Move & Number & $S_4$ Move & Number\\\hline
$\tabl{123\\231\\312} - \tabl{132\\213\\321}$ & 1 & $\tabl{1234\\2143} - \tabl{1243\\2134}$ & 18 \\
    &     &$\tabl{2314\\2431\\4123} - \tabl{2134\\2413\\4321}$ & 144\\
    &     &$\tabl{1324\\2134\\3214} - \tabl{1234\\2314\\3124}$ & 16\\
\hline
\end{tabular}
\end{table}

Although the calculation for $S_6$ is currently not possible using Gr\"obner
basis methods, there is a natural group action that reduces the
complexity of this problem. 
The group $S_n \times S_n$ acts on $\Nset^{n^2}$ by permuting rows and columns.
If we permute the rows and columns of a magic square, we still have a magic square,
therefore, this action lifts to a group action on the Markov basis of $I_{S_n}$.
In terms of tableaux, one copy of $S_n$ acts by permuting columns of the tableau,
the other acts by permuting the labels in the tableau.
We have calculated orbits under this action,
notice that the symmetrized bases are remarkably small (Table~\ref{tab:basis}).

To calculate a Markov basis for $I_{S_6}$, we had to construct the
fiber over every magic square with sum at most 5 (by
Theorem~\ref{thm:1}) and then pick moves such that every fiber is
connected by these moves \citep[see][Theorem~5.3]{gbcp}.  For degrees
2 and 3 this was relatively straightforward (e.g., there are
20,933,840 six by six magic squares with sum 3).  For these degrees,
we constructed all squares and then calculated orbits of the group
action and calculated the fiber for each orbit (there were 11 orbits
in degree 2 and 103 in degree 3).

However, there are 1,047,649,905 six by six magic squares of degree 4
and 30,767,936,616 of degree 5 \citep[from][]{ehrhart}, so complete
enumeration was not possible.  Instead, we first randomly generated
millions of magic squares with sums 4 or 5 using another Markov chain.
We broke these down into orbits, keeping track of the number of
squares we had found.  For example, we needed to generate 30 million
squares of degree 5 to find a representative for each orbit.  We were
left with 2804 orbits for degree 4 and 65481 orbits for degree
5.  For degree 5, the proof of Theorem~\ref{thm:1} shows that we only
need to consider magic squares with norm squared less that 50, leaving
13196 orbits to check.  The fibers were calculated by a depth first
search with pruning.  Remarkably, the computation showed that $I_{S_6}$ is 
generated in degree 3, see Table~\ref{tab:basis}.  

The entire calculation for $S_6$ took about 2 weeks, with the vast
majority of the time spent calculating orbits of degree 5 squares.
Our data and  code (in perl) are available for download at 
{\tt http://math.berkeley.edu/\~{}eriksson}.  The code could be easily adapted to calculate
other Markov bases with a good degree bound and a large symmetry group.
Our calculations and Table~\ref{tab:basis} suggest the following conjecture:
\begin{conj}\label{conj:1}
The ideal $I_{S_n}$ is generated in degree 3.
\end{conj}

\begin{table}
\caption{Number of generators and symmetry classes of generators by degree
in a Markov basis for $I_{S_n}$.}
\label{tab:basis}
\centering
\begin{tabular}{|c|cc|cc|cc|cc|cc|}
\hline
 &\multicolumn{2}{c|}{Degree 2} &\multicolumn{2}{c|}{Degree 3} 
&\multicolumn{2}{c|}{Degree 4} &\multicolumn{2}{c|}{Degree 5}  &\multicolumn{2}{c|}{Degree 6} \\
n & all & sym & all & sym & all & sym & all & sym & all & sym\\\hline
3 & 0 & 0 & 1 & 1&&&&&&\\
4 & 18 & 1 & 160 & 2 & 0 & 0&&&&\\
5 & 1050& 2 & 28840 & 12 & 0 & 0 & 0 & 0&&\\
6 & 57150  & 7 & 7056240 & 51 & 0 & 0 & 0 & 0 & 0 & 0\\\hline
\end{tabular}
\end{table}

\section{Structure of the toric ideal $I_{S_n}$}\label{sec:math}

Theorem 6.1 of \cite{ds} shows that every reverse lexicographic Gr\"obner basis of
$I_{S_n}$ has degree at most  $n$.  
By considering only minimal generators and not a full Gr\"obner basis, we are able
to strengthen this degree bound.

\begin{thm}\label{thm:1}
The ideal $I_{S_n}$ is generated in degree $n-1$  for $n > 3$.
\end{thm}
\begin{pf}
Since we know that $I_{S_n}$ is generated in degree $n$, we need to
show that the fibers over all magic squares with sum $n$ are each
connected by moves of degree $n-1$ or less.  Let $S$ and $T$ be tableaux in
$\rho^{-1}(\bb)$, where $\bb$ is a magic square with sum $n$.  Suppose
that the first row of $S$ and the first row of $T$ differ in exactly
$k$ places.  Then we claim that there is a degree $k+1$ move that can
be applied to $S$ to get a tableaux $S'\in \rho^{-1}(\bb)$ with the same first row as
$T$.

To change the first row of $S$ to make it agree with the first row of
$T$, we have to permute $k$ elements of the first row of $S$.  But to
remain in the fiber, this means we must also permute (at most) $k$ other
rows of $S$.  For example, if the first row of $S$ is $123\dots n$ and
the first row of $T$ is $213\dots n$, we would also have to pick the
row of $S$ with a 2 in the first column and the row with a 1 in the
second column.  Once we have picked the (at most) $k$ rows of $S$ that
must be changed, it follows from Birkhoff's
theorem \citep[e.g.,][Theorem~5.5]{vanlint} that we can change these $k$ rows and the first row to make a
new tableau $S'\in \rho^{-1}(\bb)$ that agrees with $T$ in one row.

We applied a degree $k+1$ move and are left with $S'$ and $T$ being
connected by a degree $n-1$ move, so as long as we have $k+1 \leq n -
1$, we are done.  That is, for every pair $(S,T)$ of tableaux in a
degree $n$ fiber, we must show that there is a row of $S$ and a row of
$T$ that differ in at most $n-2$ places.

Given such a pair $(S,T)$, introduce an $n \times n$ matrix $M$ where
the entries $M_{ij}$ are the number of entries that row $i$ of $S$ and
row $j$ of $T$ agree.  Notice that if $M_{ij} \geq 2$, we have rows
$i$ in $S$ and $j$ in $T$ that differ in at most $n-2$ places and are done.

Suppose that row $i$ of $S$ is $(\pi_i(1), \dots, \pi_i(n))$.  The row
sum $\sum_{j=1}^{n} M_{ij}$ counts the total number of times that
$\pi_i(j)$ appears in column $j$ for each $j$.  This is exactly
$\sum_{k=1}^{n} \bb(k, \pi(k))$.  Summing over all rows, we
see that every entry of $\bb$ gets counted its cardinality number of
times.  That is,
\[
\sum_{1\leq i,j \leq n} M_{ij} = \sum_{1\leq i,j \leq n} \bb(i,j)^2 = ||\bb||^2
\]
Now since each row of $\bb$ sums to $n$, we have that $||\bb||^2
\geq n^2$, with equality only if $\bb(i,j) = 1$ for all $i,j$.  If
this $||\bb||^2 > n^2$, then one of the $M_{ij}$ must be larger
than 1, and we are done.

Therefore, we only have to consider the fiber over
$\bb_1 = 
\left (
\begin{smallmatrix}
1 & 1 & \dots & 1 \\
\vdots & & & \vdots\\
1 & 1 & \dots & 1 \\
\end{smallmatrix}
\right ).
$
Elements of this fiber are tableaux such that every row and every
column is a permutation of $\{1, \dots, n\}$ (``Latin squares'').  Two tableaux are
connected by a degree $n-1$ move if they have a row in common.  We claim that if $n > 3$,
this graph is connected.  (Note that for $n=3$, there are two
components and a degree 3 move for $S_3$, see Table~\ref{tab:basis}.)

For fixed $\nu \in S_n$, the set $T_\nu$ of all tableaux in $\rho^{-1}(\bb_1)$
that have $\nu$ as a row is connected by definition.
Form the graph $G_n$ where the vertices are elements
$\nu \in S_n$ and there is an edge between $\lambda$ and $\nu$ if
$\lambda$ and $\nu$ occur in a tableau together.
Then if this graph is connected, the whole fiber over $\bb_1$ is connected by degree $n-1$ moves.

First, we claim that $\lambda$ and $\nu$ occur together in a tableau 
if and only if
$\lambda$ is a derangement with respect to $\nu$ (i.e., if $\lambda$ and $\nu$ 
are disjoint from each other).   The derangement condition is clearly necessary.
Sufficiency follows from Birkhoff's theorem: 
if $\lambda$ is a derangement with respect to $\nu$, then
the square $\bb_1 - \rho(\lambda) - \rho(\nu)$ has non-negative entries and row and column
sums $n-2$, therefore, it it the sum of $n-2$ permutation matrices.
Thus, $G_n$ is the graph where two permutations are connected by an edge when they
are disjoint.

Now note that 
$[1, 2, \dots, n-2, n-1, n]$ and $[3, 4, \dots, n, 1, 2]$ are connected 
in $G_n$ since the second is a cyclic shift of the first.
Then, if $n > 3$, 
$[3,4, \dots, n, 1,2]$ and  $[1,2, \dots, n-2, n, n-1]$ are also connected.
Thus $[1, 2, \dots, n]$ and $[1,2,\dots,n-2,n,n-1]$ are connected, so applying
transpositions keeps us in the same connected component of $G_n$.  But $S_n$ is 
generated by transpositions, so $G_n$ is connected and therefore $\rho^{-1}(\bb_1)$
is connected by moves of degree $n-1$.
\qed
\end{pf}

\begin{rem}
From partial computations with {\tt CaTS}
\citep{cats} for $n = 4$, it appears that every Gr\"obner basis for
$S_4$ contains degree 4 elements, while the Markov basis for $S_4$
needs only degree 3.  Furthermore, our Gr\"obner basis for $S_5$ 
contained degree 5 elements.
Therefore, it is possible that the degree $n$ Gr\"obner basis of  \cite{ds} is
the Gr\"obner basis of smallest degree.
\end{rem}

While $I_{S_n}$ is difficult to compute,
it is easy to classify the degree 2 part of the Markov basis.  
To do so, first assume that all entries of the magic square $\bb$ are either 1 or 0.
Then the squares with non-trivial $\rho^{-1}(\bb)$ are those that can
be put in a block diagonal form with $k \geq 2$ blocks and each block
of size at least 2.  Such a magic square has a fiber of size
$2^{k-1}$, corresponding to choosing, for each block,
 an orientation of the two
permutations that sum to that block  (since the order of the rows in a
tableau don't matter, there are only $k-1$ such choices).
Therefore, we need $2^{k-1} -1$ moves to make such a fiber connected.
It is a standard fact \citep[e.g.,][Chapter~1]{stanley1} that the number of partitions of $n$ into $k$ blocks
each of size at least 2 (denoted $p_2(n;k)$) satisfies
\[
\sum_{n\geq 0} p_2(n;k) q^n = 
q^{2k} \prod_{i=1}^{k} \frac 1 {1-q^i}
\]
Then let $D_2(n)$ be the number of degree 2 moves, up to symmetry, in 
a Markov basis for $S_n$.  If a magic square contains a 2, it can be thought
of as coming from $D_2(n-1)$, so putting everything together, we have
\[
D_2(n) = D_2(n-1) + \sum_{k=2}^{\lfloor \frac n 2 \rfloor} (2^{k-1}-1) [q^{n-2k}]  \prod_{i=1}^{k} \frac 1 {1-q^i},
\]
where $[q^j](\sum a_i q^i) := a_j$.
For example, $D_2(9) = 47$.  %In particular, $D_2(n) = \Omega(2^n)$.

\section{Statistical analysis of the election data}\label{sec:stat}

In order 
to run  a Markov chain fixing $\hf(\rho)$ on data $f$, we use the Markov basis 
$\{g_1, \dots, g_B\}$ as calculated above.
Then, starting from $f$, choose $i$ uniformly in $\{1,
2, \dots, B\}$ and choose $\epsilon = \pm 1$ with probability $1/2$.
If $f + \epsilon g_i \geq 0$ (coordinate-wise), the Markov chain moves
to $f + \epsilon g_i$.  Otherwise, the Markov chain stays at $f$.
This gives a symmetric connected Markov chain on the data sets with a
fixed value of $\hf(\rho)$.  As such, it has a uniform stationary
distribution.  To get a sample from the hypergeometric distribution
\eqn{eq:hyper}, the Metropolis algorithm or the Gibbs sampler can be
used \citep[see][]{liu}.

Given a symmetrized basis, we can still perform a random walk.  Pick,
at random, an element $g$ of $S_n \times S_n$.  Pick a move from the
symmetrized basis at random, apply $g$ to it (permuting columns and
renaming entries), then use the resulting move in the Markov chain.
This again gives a symmetric Markov chain that converges to the
uniform distribution.

\begin{table}
\caption{Squared length (divided by 120) of the projection of the APA data into the 7 isotypic subspaces of $S_5$.
Also, the averages of this projection for 100 random draws for 3 perturbations.}
\label{tab:sum2}
\centering
\begin{tabular}{|c|ccccccc|}
\hline
& $S^{5}$ & $S^{4,1}$ &  $S^{3,2}$ & $S^{3,1,1}$ & 
$S^{2,2,1}$ & $S^{2,1,1,1}$ & $S^{1,1,1,1,1}$\\\hline
Data      & 2286 & 298 & 459 & 78 & 27 & 7 & 0\\
Hypergeometric & 2286 & 298 & 16 & 19 & 10 & 6 & 0\\
Uniform & 2286 & 298 & 511 & 672 & 436 & 295 & 25\\
Bootstrap & 2286 & 303 & 469 &93 & 37 & 13 & 1\\
\hline
\end{tabular}
\end{table}

In this section, we apply the Markov basis for $S_5$ to analyze Table~\ref{tab:1}.
The second and third rows of Table~\ref{tab:sum2} show the average sum
of squares for 100 samples from the hypergeometric distribution
\eqn{eq:hyper} (row 2) and from the uniform distribution (row 3) with
$\hf(\rho)$ fixed.  Both sets of numbers are based on a Markov chain
simulation using a symmetrized version of the minimal basis.  In each
case, starting from the original data set, the chain was run 10,000
steps and the current function recorded.  From here, the chain was run
10,000 further steps, and so on until 100 functions were recorded.
While the running time of 10,000 steps is arbitrary, wide variation 
in the running time did not appreciably change the results.

\begin{figure}
\caption{Distribution of the length of the projection to $S^{3,2}$ with the
Metropolis and uniform random walks.}
\label{fig:met}
\centering
\includegraphics[width=2.4in]{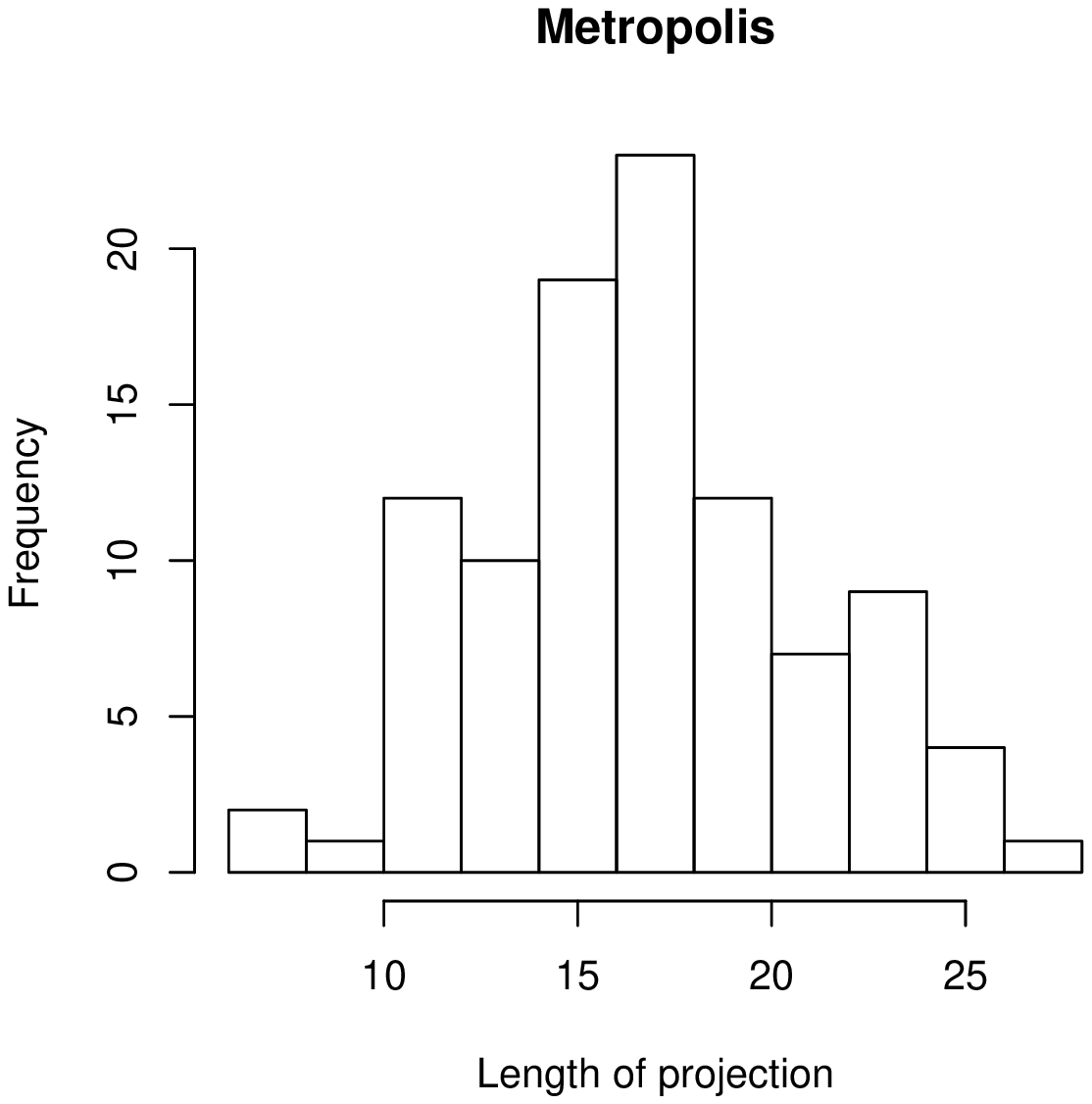}
\includegraphics[width=2.4in]{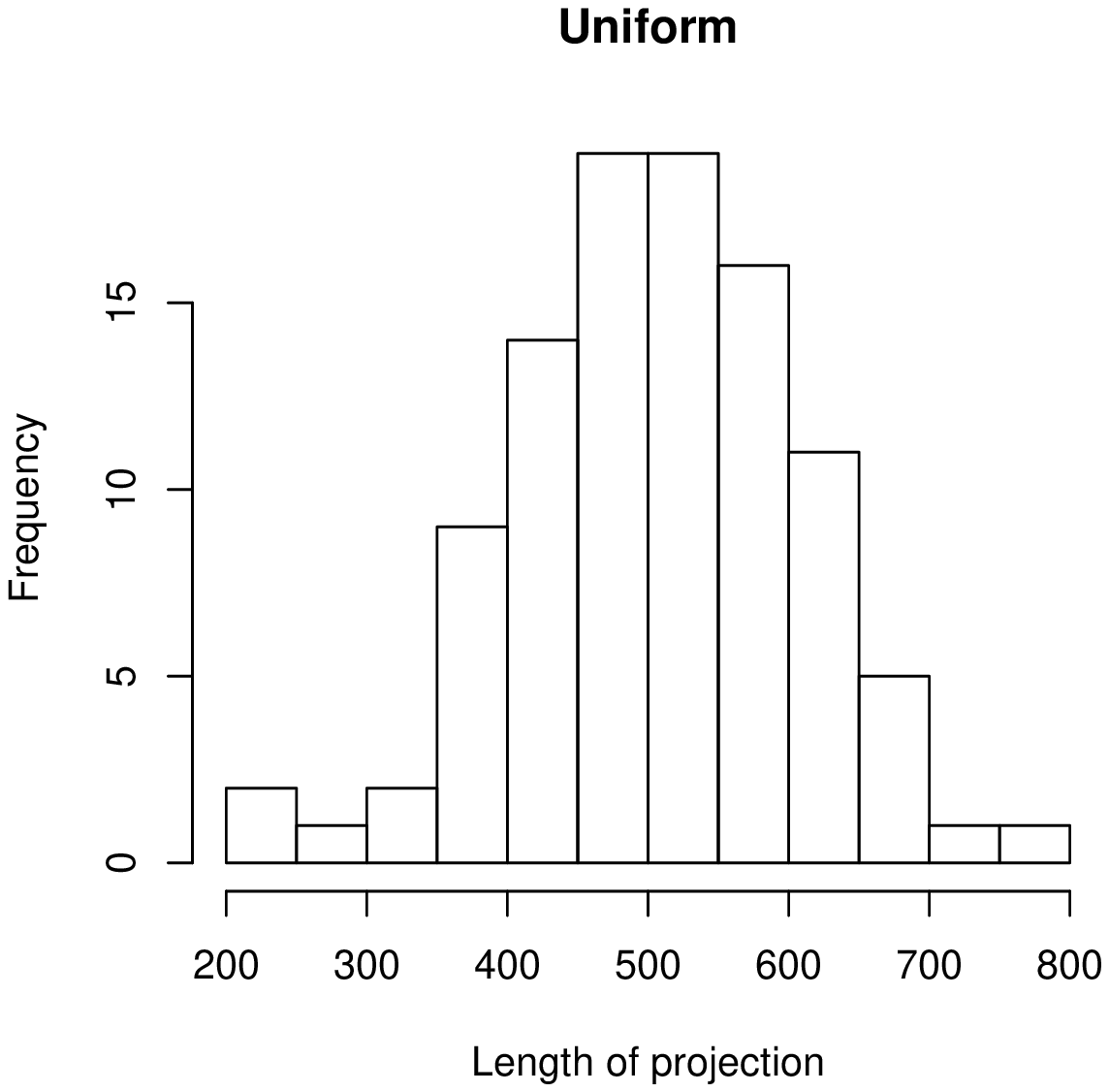}
\end{figure}

A histogram of the 100 values of the length of the projection into
$S^{3,2}$ under each distribution is shown in Figure~\ref{fig:met}.
These show some of variability but nothing exceptional.
The histograms for the other projections are very similar.

Consider first the hypergeometric distribution leading to row 2 of
Table~\ref{tab:sum2} and Figure~\ref{fig:met}.  A natural test of
goodness of fit of the model \eqn{eq:model} for the APA data may be
based on the conditional distribution of the squared length of the
projection of the data into $S^{3,2}$.  From the random walk under the
null model, this should be about $15 \pm 5$.  For the actual data,
this projection is 459.  This gives a definite reason to reject the
null model.  Our look at the data projected into $S^{3,2}$ and the
analysis that emerged in Section~\ref{sec:fourier} confirms this conclusion.

In \cite{de}, the uniform distribution of the data conditional on a
sufficient statistic was suggested as an antagonistic alternative to
the null hypothesis when the data strongly rejects a null model.  The
idea is to help quantify if the data is really far from the null, or
practically close to the null and just rejected because of a small
deviation but a large sample size \citep[see the discussion in the last
section of][]{de}.  From Figure~\ref{fig:met}, we see that the actual
projected length 459 is roughly typical of a pick from the uniform.
This affirms the strong rejection of \eqn{eq:model} and points to a
need to look at the structure of the higher order projection on its
own terms.

An appropriate stability
analysis was left open in \citet{diaconis-spectral}.  If the data in
Table~\ref{tab:1} were a sample from a larger population, the sampling
variability adds noise to the signal.  How stable is the analysis
above to natural stochastic perturbations?
One standard approach is shown in the last row of
Table~\ref{tab:sum2}.  This is based on a boot-strap perturbation of
the data in Table~\ref{tab:1}.  Here, the votes of all 5972 rankers
are put in a hat and a sample of size 5972 is drawn from the hat with
replacement to give a new data set.  The sum of squares decomposition
is repeated.  This resampling step (from the original population) was
repeated 100 times.  The entries in the last row of
Table~\ref{tab:sum2} show the average squared length of these
projections.  We see that they do not vary much from the original
sum of squares.  While not reported here, the boot-strap analogue of
the second order analysis in Table~\ref{tab:sum2} was quite stable.  We
conclude that sampling variability is not an important issue for this
example.

In \cite{ds} an $S_4$ example was analyzed.  However, the data was
analyzed using only the uniform distribution, which 
only tells half of the story. The 
analysis under hypergeometric sampling gives an important supplement.
Briefly, a sample of 2262
German citizens were asked to rank order the desirability of four
political goals.  The data and a first order summary appears in \cite{ds}.
The sizes of the projections for the data and the random walks appear in Table~\ref{tab:s4sum}.  
We have noted a typographical error in the data, the 2431 entry should be 59.

\begin{table}
\caption{Length of projections onto the 5 isotypic subspaces for the $S_4$ data and
three tests.}
\label{tab:s4sum}
\centering
\begin{tabular}{|c|ccccc|}
\hline
& $S^4$ & $S^{3,1}$ & $S^{2,2}$ &  $S^{2,1,1}$ & $S^{1,1,1,1}$\\\hline
Data      &  462 & 381 & 268 & 49 & 4 \\
Metropolis & 462 & 381 & 169 & 37 & 8 \\
Uniform &    462 & 381 & 277 & 228 & 80\\
Bootstrap &  462 & 381 & 269 & 56 & 7 \\
\hline
\end{tabular}
\end{table}

The projection of the data into the second order subspace $S^{2,2}$ has squared
length 268.  The boot-strap analysis (Line 4 in Table~\ref{tab:s4sum})
shows this is stable under sampling perturbations.  The hypergeometric
analysis (line 2 of Table~\ref{tab:s4sum}) suggests that for the
specific data, relatively large projections
onto the second order space are typical, even if the first order model holds.  
This is quite different than the previous example.  Still,
the observed 268 is sufficiently much larger than 169 
that a look at the second order projection is warranted.
The uniform analysis points to the actual projection being typical,
this again suggests a serious look at the second order projection.

As a side remark, the software {\tt LattE} \citep{latte} 
can be used to count how many data sets have a given
first order summary.  For our $S_4$ example, these
correspond to lattice points inside a 
convex polytope with 6285 vertices in $\Rset^{24}$.  
{\tt LattE} computes (in only 523.12 seconds) that there are
11606690287805167142987310121 (approximately 
$10^{28}$) elements of $\Nset[S_4]$ with 
the same first order summary as our $S_4$ example.

\section{Conclusions}
In this paper, we have given a general methodology for studying group
valued data where the summary we are interested in is given by a representation 
of the group and analyzed in detail the case of ranked data.  
This suggests a family of interesting toric ideals: to each
finite group $G$ and representation $\rho$ we associate a
toric ideal  (Definition~\ref{def:ideal}).

For practical purposes, it would be nice to have a general algorithm
to analyze ranked data with $n$ candidates.  We ran Markov chains
using just the degree 2 moves, but they seemed to mix very poorly.
However, our computations and Conjecture~\ref{conj:1} suggest that
finding all (or even some) degree 3 moves in addition to the degree 2
moves would allow for a good random walk.

\section*{Acknowledgments}
We thank Susan Holmes and Aaron Staple for writing the {\tt R} code that
made this analysis possible, Ruriko Yoshida for computational help with
{\tt LattE}, and Bernd Sturmfels and Seth Sullivant for helpful comments.  This
work was a direct result of an AIM workshop 
in December 2003.  We thank the organizers  and
the staff at AIM for a great workshop.

\bibliographystyle{elsart-harv}
%\bibliography{bib}

\begin{thebibliography}{18}
\expandafter\ifx\csname natexlab\endcsname\relax\def\natexlab#1{#1}\fi
\expandafter\ifx\csname url\endcsname\relax
  \def\url#1{\texttt{#1}}\fi
\expandafter\ifx\csname urlprefix\endcsname\relax\def\urlprefix{URL }\fi

\bibitem[{Beck and Pixton(2003)}]{ehrhart}
Beck, M., Pixton, D., 2003. The {E}hrhart polynomial of the {B}irkhoff
  polytope. Discrete Comput. Geom. 30~(4), 623--637.

\bibitem[{De~Loera et~al.(2003)De~Loera, Haws, Hemmecke, Huggins, Tauzer, and
  Yoshida}]{latte}
De~Loera, J., Haws, D., Hemmecke, R., Huggins, P., Tauzer, J., Yoshida, R.,
  2003. A user's guide for latte v1.1. Available at {\tt
  http://www.math.ucdavis.edu/\~{}latte/}.

\bibitem[{Diaconis(1988)}]{diaconis-book}
Diaconis, P., 1988. Group Representations in Probability and Statistics.
  Vol.~11 of IMS Lecture Series. Institute of Mathematical Statistics.

\bibitem[{Diaconis(1989)}]{diaconis-spectral}
Diaconis, P., 1989. A generalization of spectral analysis with application to
  ranked data. Ann. Statist. 17~(3), 949--979.

\bibitem[{Diaconis and Efron(1985)}]{de}
Diaconis, P., Efron, B., 1985. Testing for independence in a two-way table: new
  interpretations of the chi-square statistic. Ann. Statist. 13~(3), 845--913.

\bibitem[{Diaconis and Freedman(1984)}]{diaconis-freedman}
Diaconis, P., Freedman, D., 1984. Partial exchangeability and sufficiency. In:
  Ghosh, J., Roy, J. (Eds.), Statistics: Applications and New Directions.
  Indian Statistical Institute, Calcutta, pp. 205--236.

\bibitem[{Diaconis and Sturmfels(1998)}]{ds}
Diaconis, P., Sturmfels, B., 1998. Algebraic algorithms for sampling from
  conditional distributions. Ann. Statist. 26~(1), 363--397.

\bibitem[{Fisher(1973)}]{fisher}
Fisher, R.~A., 1973. Statistical methods for research workers. Hafner
  Publishing Co., New York.

\bibitem[{Hemmecke and Hemmecke(2003)}]{4ti2}
Hemmecke, R., Hemmecke, R., Sep. 2003. 4ti2 version 1.1---computation of
  {H}ilbert bases, {G}raver bases, toric {G}r{\"o}bner bases, and more.
  Available at www.4ti2.de.

\bibitem[{Jensen(2003)}]{cats}
Jensen, A.~N., 2003. Cats, a software system for toric state polytopes.
  Available at {\tt http://www.soopadoopa.dk/anders/cats/cats.html}.

\bibitem[{Liu(2001)}]{liu}
Liu, J., 2001. Monte Carlo techniques in scientific computing. Springer, New
  York.

\bibitem[{Marden(1995)}]{marden}
Marden, J., 1995. Analyzing and modeling rank data. Chapman and Hall, London.

\bibitem[{Silverberg(1984)}]{silver}
Silverberg, A., 1984. Statistical models for $q$-permutations. Proc Biopharm.
  Sec. Amer. Statist. Assoc., 107--112.

\bibitem[{Stanley(1997)}]{stanley1}
Stanley, R.~P., 1997. Enumerative combinatorics. {V}ol. 1. Vol.~49 of Cambridge
  Studies in Advanced Mathematics. Cambridge University Press, Cambridge.

\bibitem[{Stern(1993)}]{stern}
Stern, H., 1993. Probability models on rankings and the electoral process. In:
  Fligner, M., Verducci, J. (Eds.), Probability models and statistical analyses
  for ranking data. Springer, pp. 173--195.

\bibitem[{Sturmfels(1996)}]{gbcp}
Sturmfels, B., 1996. Gr\"obner bases and convex polytopes. Vol.~8 of University
  Lecture Series. American Mathematical Society, Providence, RI.

\bibitem[{van Lint and Wilson(2001)}]{vanlint}
van Lint, J.~H., Wilson, R.~M., 2001. A course in combinatorics, 2nd Edition.
  Cambridge University Press, Cambridge.

\bibitem[{Verducci(1982)}]{verd}
Verducci, J., 1982. Discriminating between two probabilities on the basis of
  ranked preferences. Ph.D. thesis, Stanford University.

\end{thebibliography}

\end{document}